\documentclass{amsart}
\usepackage{amsmath,amssymb,graphicx,latexsym}
\usepackage{hyperref}
\usepackage[T1]{fontenc}
\usepackage{textcomp}
\usepackage[utf8]{inputenc}
\usepackage{mathtools} 
\usepackage{amsfonts} 
\usepackage{times, url}
\usepackage{siunitx}
\usepackage{float}
\usepackage{listings}
\usepackage{subcaption}
\usepackage{nccmath}
\usepackage[english]{babel}

\theoremstyle{definition}

\newtheorem*{definition}{Definition}
\newtheorem*{notation*}{Notation}

\begin{document}
\title{Observations regarding the repetition of the last digits of a tetration of generic base}
\author{Luca Onnis}
\date{September 2021}
\maketitle
\begin{abstract}
    This paper investigates the behavior of the last digits of a tetration \cite{tet}of generic base. In fact, last digits of a tetration are the same starting from a certain hyper-exponent and in order to compute them we reduce those expressions $\mod 10^{n}$. Very surprisingly (although unproved) I think that the repeating digits depend on the residue $\mod 10$ of the base and on the exponents of a particular way to express that base. Then I'll discuss about the results and I'll show different tables and examples in order to support my conjecture.
\end{abstract}
\section{Definitions}
\begin{definition}[Tetration]
In mathematics, tetration is an operation based on iterated, or repeated, exponentiation. It is the next hyperoperation after exponentiation, but before pentation. The word was coined by Reuben Louis Goodstein from tetra- (four) and iteration. ${^{a}n}$ represent the $a$-th tetration of $n$ , or:
\[
n^{n^{n^{n^{\dots}}}}\Bigl\} \mbox{ $a$} \mbox{ times}
\]
\end{definition}
\begin{definition}[Floor and Ceiling function] \cite{tdn}
In mathematics and computer science, the floor function is the function that takes as input a real number x, and gives as output the greatest integer less than or equal to x, denoted $\lfloor x \rfloor$. Similarly, the ceiling function maps x to the least integer greater than or equal to x, denoted $\lceil x \rceil$.
\end{definition}
\section{Main conjecture}
Let $f_{q}(x,y,n)$ be a function such that if:
\[
f_{q}(x,y,n)=u
\]
Then:
\[
{^{\infty}\Bigl[q^{(2^{x}\cdot5^{y})\cdot a}\Bigr]} \equiv {^{u}\Bigl[q^{(2^{x}\cdot5^{y})\cdot a}\Bigr]} \mod (10^{n})
\]
where $x,y,n,q,a \in\mathbb{N}$ , $q\not=10h$, $a \not=2h$ and $a\not=5h$ and $u$ is the minimum value such that this congruence is true. \\
\textbf{Note that these formulas work if $x\geq2$} \\
At the end of this paper there will be a section for the $x<2$ case.
\subsection{$q\equiv 1,9 \mod 10$}
I define $\Delta_2$ and $\Delta_5$ as:
\[
\Delta_2=\max[v_2(q+1),v_2(q-1)]
\]
\[
\Delta_5=\max[v_5(q+1),v_5(q-1)]
\]
We'll have that:
\[
f_{q \equiv 1,9 \mod 10}(x,y,n)=\max\Biggl[\Bigl\lceil\frac{n}{x+\Delta_2}\Bigr\rceil,\Bigl\lceil\frac{n}{y+\Delta_5}\Bigr\rceil\Biggr]-1
\]
\subsection{$q\equiv 3,7 \mod 10$}
I define $\Gamma_2$ and $\Gamma_5$ as:
\[
\Gamma_2=\max[v_2(q+1),v_2(q-1)]
\]
\[
\Gamma_5=\max[v_5(q^{2}+1),v_5(q^{2}-1)]
\]
We'll have that:
\[
f_{q \equiv 3,7 \mod 10}(x,y,n)=\max\Biggl[\Bigl\lceil\frac{n}{x+\Gamma_2}\Bigr\rceil,\Bigl\lceil\frac{n}{y+\Gamma_5}\Bigr\rceil\Biggr]-1
\]
\subsection{$q\equiv 5 \mod 10$}
$q=10h+5$ for some $h\in\mathbb{N}$. It is also possible to write it as $q=5^{k}h$ for some $h,k\in\mathbb{N}$. As we know from my previous paper, if the base is a power of 5, the last repeating digits only depend on the 2-adic valuation of the base's exponent. This is right not only for the power of base 5 but for all multiples of 5. We can write the base as $5^{2^{x}k}h$, for some $h,k\in\mathbb{N}$ , $k\not=2a\wedge h\not=5b$. So we'll have:
\[
f_{q \equiv 5 \mod 10}(x,y,n)= \Bigl\lceil\frac{n}{x+\Delta_2}\Bigr\rceil-1
\]
Note that here $y$ could be every integer number, but the result depends on $x$.
\subsection{$q\equiv 0 \mod 2$}
$q=10h+2$ for some $h\in\mathbb{N}$. It is also possible to write it as $q=2^{k}h$ for some $h,k\in\mathbb{N}$. As we know from my previous paper, if the base is a power of 2, the last repeating digits only depend on the 5-adic valuation of the base's exponent. This is right not only for the power of base 2 but for all multiples of 2. We can write the base as $2^{5^{y}k}h$, for some $h,k\in\mathbb{N}$ , $k\not=5a\wedge h\not=2b$. So we'll have:
\[
f_{q \equiv 0 \mod 2}(x,y,n)= \Bigl\lceil\frac{n}{y+\Gamma_5}\Bigr\rceil-1
\]
Note that here $x$ could be every integer number, but the result depends on $y$.
\section{Introduction}
When you hear talking about tetrations , you instantly think about a giant number, which is very difficult to compute. One of the first things you could notice when you are studying them is that the last digits of every tetration with a positive integer base begins to repeat after a certain hyper-exponent.
\section{Examples and other observations}
First of all I would like to mention that the $q\equiv 1,3,7,9 \mod 10$ case is now divided in 2 sub-cases, namely: $q\equiv 1,9 \mod 10$ and $q\equiv 3,7 \mod 10$. Then I understood that they can not be grouped together. Furthermore, initially I thought that the last repeating digits were depending on the prime factorization of the base $q$, but then I understood that they only depend on the residue $\mod 10$ of that base and on the relationship of the discriminants $\Delta$ and $\Gamma$. Anyway at the end of this paper it is possible to find a table of some values of $f_p(x,y,n)$ where $p$ is a prime number.
\subsection{Example 1}
Consider the infinite tetration of $4599^{2^{8}\cdot 5^{5}}$ , or $4599^{800000}$. We know from our first formula that the last 40 digits are the same starting from the 5-th tetration of that number. Indeed, $4599 \equiv 9 \mod 10$ and $\lceil\frac{n}{y+\Delta_5}\rceil\geq\lceil\frac{n}{x+\Delta_2}\rceil$. \\
In fact:
\[
\Delta_2=\max[v_2(4599+1),v_2(4599-1)]=\max[3,1]=3
\]
\[
\Delta_5=\max[v_5(4599+1),v_5(4599-1)]=\max[2,0]=2
\]
And:
\[
\Bigl\lceil\frac{40}{5+2}\Bigr\rceil\geq\Bigl\lceil\frac{40}{8+3}\Bigr\rceil
\]
So we'll have that:
\[
f_{4599}(8,5,40)=\Bigl\lceil\frac{40}{5+\max[v_5(4600),v_5(4598)]}\Bigr\rceil-1
\]
\[
f_{4599}(8,5,40)=\Bigl\lceil\frac{40}{5+\max[2,1]}\Bigr\rceil-1
\]
\[
f_{4599}(8,5,40)=\Bigl\lceil\frac{40}{5+2}\Bigr\rceil-1
\]
\[
f_{4599}(8,5,40)=5
\]
In fact:
\[
 {^{5}\Bigl[4599^{(2^{8}\cdot5^{5})}\Bigr]} \equiv 574081590929428693334403581932320000001 \mod (10^{40})
\]
And also:
\[
 {^{6}\Bigl[4599^{(2^{8}\cdot5^{5})}\Bigr]} \equiv 574081590929428693334403581932320000001 \mod (10^{40})
\]
And so on for every hyper-exponent greater or equal to 5. So 5 is minimum, in fact if you consider the fourth tetration of that base, it doesn't work.
\[
 {^{4}\Bigl[4599^{(2^{8}\cdot5^{5})}\Bigr]} \equiv 3530881590929428693334403581932320000001 \mod (10^{40})
\]
\subsection{Example 2}
Consider the infinite tetration of $1251^{2^{2}\cdot 5^{4}}$ , or $1251^{2500}$. We know from our first formula that the last 30 digits are the same starting from the 7-th tetration of that number. Indeed, $1251 \equiv 1 \mod 10$ and $\lceil\frac{n}{x+\Delta_2}\rceil\geq\lceil\frac{n}{y+\Delta_5}\rceil$. \\
In fact:
\[
\Delta_2=\max[v_2(1251+1),v_2(1251-1)]=\max[2,1]=2
\]
\[
\Delta_5=\max[v_5(1251+1),v_5(1251-1)]=\max[0,4]=4
\]
And:
\[
\Bigl\lceil\frac{30}{2+2}\Bigr\rceil\geq\Bigl\lceil\frac{30}{4+4}\Bigr\rceil
\]
So we'll have that:
\[
f_{1251}(2,4,30)=\Bigl\lceil\frac{30}{2+\max[v_2(1252),v_2(1250)]}\Bigr\rceil-1
\]
\[
f_{1251}(2,4,30)=\Bigl\lceil\frac{30}{2+2}\Bigr\rceil-1=7
\]
In fact:
\[
 {^{7}\Bigl[1251^{(2^{2}\cdot5^{4})}\Bigr]} \equiv 297934155568039465330081250001 \mod (10^{30})
\]
And also:
\[
 {^{8}\Bigl[1251^{(2^{2}\cdot5^{4})}\Bigr]} \equiv 297934155568039465330081250001 \mod (10^{30})
\]
And so on for every hyper-exponent greater or equal to 7. So 7 is minimum, in fact if you consider the sixth tetration of that base, it doesn't work.
\[
 {^{6}\Bigl[1251^{(2^{2}\cdot5^{4})}\Bigr]} \equiv 47934155568039465330081250001 \mod (10^{30})
\]
\subsection{Example 3}
Consider the infinite tetration of $17^{2^{3}\cdot 5^{6}}$ , or $17^{125000}$. We know from our second formula that the last 53 digits are the same starting from the 7-th tetration of that number. Indeed, $17 \equiv 7 \mod 10$ and $\lceil\frac{n}{y+\Gamma_5}\rceil\geq\lceil\frac{n}{x+\Gamma_2}\rceil$. \\
In fact:
\[
\Gamma_2=\max[v_2(17+1),v_2(17-1)]=\max[1,4]=4
\]
\[
\Gamma_5=v_5(17^{2}+1)=1
\]
And:
\[
\Bigr\lceil\frac{53}{6+1}\Bigr\rceil\geq\Bigl\lceil\frac{53}{3+4}\Bigr\rceil
\]
So we'll have that:
\[
f_{17}(3,6,53)=\Bigl\lceil\frac{53}{6+v_5(290)}\Bigl\rceil-1
\]
\[
f_{17}(3,6,53)=\Bigl\lceil\frac{53}{7}\Bigl\rceil-1=7
\]
So we'll have that:
\[
 {^{7}\Bigl[17^{(2^{3}\cdot5^{6})}\Bigr]} \equiv 52737008157199929548933683973150858896289457010000001 \mod (10^{53})
\]
And so on for every hyper-exponent greater or equal to 7.
\subsection{Example 4}
Consider the infinite tetration of $63^{2^{5}\cdot 5^{2}\cdot 3}$ , or $63^{2400}$. We know from our second formula that the last 15 digits are the same starting from the 4-th tetration of that number. Indeed, $63 \equiv 3 \mod 10$ and $\lceil\frac{n}{y+\Gamma_5}\rceil\geq\lceil\frac{n}{x+\Gamma_2}\rceil$. \\
In fact:
\[
\Gamma_2=\max[v_2(63+1),v_2(63-1)]=\max[6,1]=6
\]
\[
\Gamma_5=v_5(63^{2}+1)=1
\]
And:
\[
\Bigr\lceil\frac{15}{2+1}\Bigr\rceil\geq\Bigl\lceil\frac{15}{5+6}\Bigr\rceil
\]
So we'll have that:
\[
f_{63}(5,2,15)=\Bigl\lceil\frac{15}{2+v_5(3970)}\Bigl\rceil-1
\]
\[
f_{63}(5,2,15)=\Bigl\lceil\frac{15}{3}\Bigl\rceil-1=4
\]
So we'll have that:
\[
 {^{4}\Bigl[63^{(2^{5}\cdot5^{2}\cdot 3)}\Bigr]} \equiv 547909642496001 \mod (10^{15})
\]
And so on for every hyper-exponent greater or equal to 4.
\subsection{Example 5}
Consider the infinite tetration of $255^{2^{4}\cdot 5^{3}}$ , or $255^{2000}$. We know from our third formula that the last 34 digits are the same starting from the 3-th tetration of that number. Indeed, $255 \equiv 5 \mod 10$. \\
So we'll have that:
\[
f_{255}(4,3,34)=\Bigl\lceil\frac{34}{3+\max[v_2(256),v_2(254)]}\Bigl\rceil-1
\]
\[
f_{255}(4,3,34)=\Bigl\lceil\frac{34}{11}\Bigl\rceil-1=3
\]
So:
\[
 {^{3}\Bigl[255^{(2^{4}\cdot5^{3})}\Bigr]} \equiv 6154363253735937178134918212890625 \mod (10^{34})
\]
And so on for every hyper-exponent greater or equal to 3.
\subsection{Example 6}
Consider the infinite tetration of $192^{2^{2}\cdot 5^{3}}$ , or $192^{500}$. We know from our fourth formula that the last 20 digits are the same starting from the 4-th tetration of that number. Indeed, $192 \equiv 0 \mod 2$ and:
\[
\Gamma_5=\max[v_5(192^{2}+1),v_5(192^{2}-1)]=\max[1,0]=1
\]
So we'll have that:
\[
f_{192}(2,3,20)=\Bigl\lceil\frac{20}{3+1}\Bigl\rceil-1
\]
\[
f_{192}(2,3,20)=\Bigl\lceil\frac{20}{4}\Bigl\rceil-1=4
\]
So:
\[
 {^{4}\Bigl[255^{(2^{2}\cdot5^{3})}\Bigr]} \equiv 14517958004101349376 \mod (10^{20})
\]
And so on for every hyper-exponent greater or equal to 4.
\section{$x<2$ case} 
I think that our formula in order to work must respect the condition: $x\geq2$. In fact if $x=0$ or $x=1$ , some of the results seem to be different from the main formula.
\subsection{When $q\equiv 1,9 \mod 10$}
We are now considering the case where the base $q$ is a number congruent to $1,9\mod 10$. Consider a function $f_q(n)$ such that if:
\[
f_{q\equiv 1,9 \mod 10}(x,0,n)=u
\]
Then:
\[
{^{\infty}q} \equiv {^{u}\Bigl[q^{2^{x}}\Bigr]} \mod (10^{n})
\]
\begin{table}[H]
\centering
\begin{tabular}{|| l | l | l || l | l | l ||}
 \hline
$q\equiv 1 \mod 10$ & $f_q(0,0,n)$ & $f_q(1,0,n)$ & $q\equiv 9 \mod 10$ & $f_q(0,0,n)$ & $f_q(1,0,n)$ \\
 \hline
 1 & & & 9 & $n$ & $n-1$ \\
 \hline
 11 & $n-1$ & $n-1$ & 19 & $n$ & $n-1$ \\
 \hline
 21 & $n-1$ & $n-1$ & 29 & $n$ & $n-1$ \\
 \hline
 31 & $n-1$ & $n-1$ & 39 & $n$ & $n-1$ \\
 \hline
 41 & $n-1$ & $n-1$ & 49 & $\lceil\frac{n}{2}\rceil$ & $\lceil\frac{n}{2}\rceil-1$ \\
 \hline
 51 & $\lfloor\frac{n}{2}\rfloor$ & $\lceil\frac{n}{2}\rceil-1$ & 59 & $n$ & $n-1$ \\
 \hline
 61 & $n-1$ & $n-1$ & 69 & $n$ & $n-1$ \\
 \hline
 71 & $n-1$ & $n-1$ & 79 & $n$ & $n-1$ \\
 \hline
 81 & $n-1$ & $n-1$ & 89 & $n$ & $n-1$ \\
 \hline
 91 & $n-1$ & $n-1$ & 99 & $\lceil\frac{n}{2}\rceil$ & $\lceil\frac{n}{2}\rceil-1$ \\
 \hline
 101 & $\lceil\frac{n}{2}\rceil-1$ & $\lceil\frac{n}{2}\rceil-1$ & 109 & $n$ & $n-1$ \\
 \hline
 \end{tabular}
\caption{Values of $f_q(x,0,n)$ by varying $q$ and $x$}
\label{tab:4.1}
\end{table} 
We see something particular when $q\equiv1\mod50$ and when $q\equiv49\mod50$:
\begin{table}[H]
\centering
\begin{tabular}{|| l | l | l || l | l | l ||}
 \hline
 $q\equiv 1 \mod 50$ & $f_q(0,0,n)$ & $f_q(1,0,n)$ & $q\equiv 49 \mod 50$ & $f_q(0,0,n)$ & $f_q(1,0,n)$ \\
 \hline
 1 & & & 49 & $\lceil\frac{n}{2}\rceil$ & $\lceil\frac{n}{2}\rceil-1$ \\
 \hline
 51 & $\lfloor\frac{n}{2}\rfloor$ & $\lceil\frac{n}{2}\rceil-1$ & 99 & $\lceil\frac{n}{2}\rceil$ & $\lceil\frac{n}{2}\rceil-1$ \\
 \hline
 101 & $\lceil\frac{n}{2}\rceil-1$ & $\lceil\frac{n}{2}\rceil-1$ & 149 & $\lceil\frac{n}{2}\rceil$ & $\lceil\frac{n}{2}\rceil-1$ \\
 \hline
 151 & $\lceil\frac{n}{2}\rceil-1$ & $\lceil\frac{n}{2}\rceil-1$ & 199 & $\lceil\frac{n}{2}\rceil$ & $\lceil\frac{n}{2}\rceil-1$ \\
 \hline
 201 & $\lceil\frac{n}{2}\rceil-1$ & $\lceil\frac{n}{2}\rceil-1$ & 249 & $\lceil\frac{n}{3}\rceil$ & $\lceil\frac{n}{3}\rceil-1$ \\
 \hline
 251 & $\lfloor\frac{n}{2}\rfloor$ & $\lceil\frac{n}{2}\rceil-1$ & 299 & $\lceil\frac{n}{2}\rceil$ & $\lceil\frac{n}{2}\rceil-1$ \\
 \hline
 
  \end{tabular}
\caption{Values of $f_q(x,0,n)$ by varying $q$ and $x$}
\label{tab:4.1}
\end{table} 
Using our formula we see that the last $-1$ term is fundamental, until $x\geq 1$. But if $x,y=0$, as in the first column of the table, something strange happens. \\ 
When $q\equiv 9 \mod 10$:
\[
f_{q \equiv 9 \mod 10}(0,0,n)=\max\Biggl[\Bigl\lceil\frac{n}{\Delta_2}\Bigr\rceil,\Bigl\lceil\frac{n}{\Delta_5}\Bigr\rceil\Biggr]-1
\]
When $q\equiv 1 \mod 10$:
\[
f_{q \equiv 1 \mod 10}(0,0,n)= \begin{cases} \max\Bigl[\lfloor\frac{n}{\Delta_2}\rfloor,\lfloor\frac{n}{\Delta_5}\rfloor\Bigr] & \mbox{ if} \mbox{ $q\equiv 51 \mod 200$} \\ \max\Bigl[\lceil\frac{n}{\Delta_2}\rceil,\lceil\frac{n}{\Delta_5}\rceil\Bigr]-1 & \mbox{ otherwise} \end{cases}
\]
\subsection{When $q\equiv 3,7\mod 10$}
We are now considering the case where the base $q$ is a number congruent to $3,7\mod 10$. Consider a function $f_q(n)$ such that if:
\[
f_{q\equiv 1,9 \mod 10}(x,0,n)=u
\]
Then:
\[
{^{\infty}q} \equiv {^{u}\Bigl[q^{2^{x}}\Bigr]} \mod (10^{n})
\]
\begin{table}[H]
\centering
\begin{tabular}{|| l | l | l || l | l | l ||}
 \hline
$q\equiv 3 \mod 10$ & $f_q(0,0,n)$ & $f_q(1,0,n)$ & $q\equiv 7 \mod 10$ & $f_q(0,0,n)$ & $f_q(1,0,n)$ \\
 \hline
 3 & $n+1$ & $n$ & 7 & $\lceil\frac{n}{2}\rceil+1$ & $\lceil\frac{n}{2}\rceil$ \\
 \hline
 13 & $n$ & $n$ & 17 & $n$ & $n$ \\
 \hline
 23 & $n+1$ & $n$ & 27 & $n+1$ & $n$ \\
 \hline
 33 & $n$ & $n$ & 37 & $n$ & $n$ \\
 \hline
 43 & $\lceil\frac{n}{2}\rceil+1$ & $\lceil\frac{n}{2}\rceil$ & 47 & $n+1$ & $n$ \\
 \hline
 53 & $n$ & $n$ & 57 & $\lceil\frac{n}{3}\rceil$ & $\lceil\frac{n}{3}\rceil$ \\
 \hline
 63 & $n+1$ & $n$ & 67 & $n+1$ & $n$ \\
 \hline
 73 & $n$ & $n$ & 77 & $n$ & $n$ \\
 \hline
 83 & $n+1$ & $n$ & 87 & $n+1$ & $n$ \\
 \hline
 93 & $\lceil\frac{n}{2}\rceil$ & $\lceil\frac{n}{2}\rceil$ & 97 & $n$ & $n$ \\
 \hline
 103 & $n+1$ & $n$ & 107 & $\lceil\frac{n}{2}\rceil+1$ & $\lceil\frac{n}{2}\rceil$ \\
 \hline
 \end{tabular}
\caption{Values of $f_q(x,0,n)$ by varying $q$ and $x$}
\label{tab:4.1}
\end{table} 
We see something particular when $q\equiv43\mod50$ and when $q\equiv7\mod50$:
\begin{table}[H]
\centering
\begin{tabular}{|| l | l | l || l | l | l ||}
 \hline
 $q\equiv 43 \mod 50$ & $f_q(0,0,n)$ & $f_q(1,0,n)$ & $q\equiv 7 \mod 50$ & $f_q(0,0,n)$ & $f_q(1,0,n)$ \\
 \hline
 43 & $\lceil\frac{n}{2}\rceil+1$ & $\lceil\frac{n}{2}\rceil$ & 7 & $\lceil\frac{n}{2}\rceil+1$ & $\lceil\frac{n}{2}\rceil$ \\
 \hline
 93 & $\lceil\frac{n}{2}\rceil$ & $\lceil\frac{n}{2}\rceil$ & 57 & $\lceil\frac{n}{3}\rceil$ & $\lceil\frac{n}{3}\rceil$ \\
 \hline
 143 & $\lceil\frac{n}{2}\rceil+1$ & $\lceil\frac{n}{2}\rceil$ & 107 & $\lceil\frac{n}{2}\rceil+1$ & $\lceil\frac{n}{2}\rceil$ \\
 \hline
 193 & $\lceil\frac{n}{3}\rceil$ & $\lceil\frac{n}{3}\rceil$ & 157 & $\lceil\frac{n}{2}\rceil$ & $\lceil\frac{n}{2}\rceil$ \\
 \hline
 243 & $\lceil\frac{n}{2}\rceil+1$ & $\lceil\frac{n}{2}\rceil$ & 207 & $\lceil\frac{n}{2}\rceil+1$ & $\lceil\frac{n}{2}\rceil$ \\
 \hline
 293 & $\lceil\frac{n}{2}\rceil$ & $\lceil\frac{n}{2}\rceil$ & 257 & $\lceil\frac{n}{2}\rceil$ & $\lceil\frac{n}{2}\rceil$ \\
 \hline
 343 & $\lceil\frac{n}{2}\rceil+1$ & $\lceil\frac{n}{2}\rceil$ & 307 & $\lfloor\frac{n}{2}\rfloor$ & $\lceil\frac{n}{3}\rceil$ \\
 \hline
  \end{tabular}
\caption{Values of $f_q(x,0,n)$ by varying $q$ and $x$}
\label{tab:4.1}
\end{table} 
We see something particular when $q\equiv193\mod250$ and when $q\equiv57\mod250$. This remind us some kind of $p$-adic valuation structure.
\begin{table}[H]
\centering
\begin{tabular}{|| l | l | l || l | l | l ||}
 \hline
 $q\equiv 193 \mod 250$ & $f_q(0,0,n)$ & $f_q(1,0,n)$ & $q\equiv 57 \mod 250$ & $f_q(0,0,n)$ & $f_q(1,0,n)$ \\
 \hline
 193 & $\lceil\frac{n}{3}\rceil$ & $\lceil\frac{n}{3}\rceil$ & 57 & $\lceil\frac{n}{3}\rceil$ & $\lceil\frac{n}{3}\rceil$ \\
 \hline
 443 & $\lfloor\frac{n}{2}\rfloor$ & $\lceil\frac{n}{3}\rceil-1$ & 307 & $\lfloor\frac{n}{2}\rfloor$ & $\lceil\frac{n}{3}\rceil$ \\
 \hline
  \end{tabular}
\caption{Values of $f_q(x,0,n)$ by varying $q$ and $x$}
\label{tab:4.1}
\end{table} 
It is like terms are alternating; a term with $+1$ and the other without the $+1$, looking at $f_q(0,0,n)$. Even though there are some exceptions; for example , very surprisingly, $f_{307}(0,0,n)=\lfloor\frac{n}{2}\rfloor$. 
\subsection{When $q\equiv 5\mod 10$}
When $q\equiv 5\mod 10$ our third formula seems to work fine for all integer $x\geq1$. For $x=0$, in some terms don't appear the $-1$ part of the main formula.
\subsection{Prime numbers}
Initially I was thinking that the formulas had to depend on the prime factorization of the base of the tetration, but then I have understood that our formulas only depend on the residue $\mod 10$ of the base $q$.
Here below is a table for the first 60 prime numbers:
\begin{table}[H]
\centering
\begin{tabular}{|| l | l | l | l || l | l | l | l ||}
 \hline
$p$ & $f_p(0,0,n)$ & $f_p(1,0,n)$ & $f_p(x,0,n)$ & $p$ & $f_p(0,0,n)$ & $f_p(1,0,n)$ & $f_p(x,0,n)$ \\
\hline
2 & $n+2$ & $n+1$ & $n-1$ & 127 & $n+1$ & $n$ & $n-1$ \\
 \hline
3 & $n+1$ & $n$ & $n-1$ & 131 & $n-1$ & $n-1$ & $n-1$\\
 \hline
5 & $\lceil\frac{n}{2}\rceil-1$ & $\lceil\frac{n}{3}\rceil-1$ & $\lceil\frac{n}{x+2}\rceil-1$ & 137 & $n$ & $n$ & $n-1$ \\
 \hline
7 & $\lceil\frac{n}{2}\rceil+1$ & $\lceil\frac{n}{2}\rceil$ & $\lceil\frac{n}{2}\rceil-1$ & 139 & $n$ & $n-1$ & $n-1$ \\
 \hline
11 & $n-1$ & $n-1$ & $n-1$ & 149 & $\lceil\frac{n}{2}\rceil$ & $\lceil\frac{n}{2}\rceil-1$ & $\lceil\frac{n}{2}\rceil-1$\\
 \hline
13 & $n$ & $n$ & $n-1$ & 151 & $\lceil\frac{n}{2}\rceil-1$ & $\lceil\frac{n}{2}\rceil-1$ & $\lceil\frac{n}{2}\rceil-1$\\
 \hline
17 & $n$ & $n$ & $n-1$ & 157 & $\lceil\frac{n}{2}\rceil$ & $\lceil\frac{n}{2}\rceil$ & $\lceil\frac{n}{2}\rceil-1$\\
 \hline
19 & $n$ & $n-1$ & $n-1$ & 163 & $n+1$ & $n$ & $n-1$\\
 \hline
23 & $n+1$ & $n$ & $n-1$ & 167 & $n+1$ & $n$ & $n-1$\\
 \hline
29 & $n$ & $n-1$ & $n-1$ & 173 & $n$ & $n$ & $n-1$\\
 \hline
31 & $n-1$ & $n-1$ & $n-1$ & 179 & $n$ & $n-1$ & $n-1$\\
 \hline
37 & $n$ & $n$ & $n-1$ & 181 & $n-1$ & $n-1$ & $n-1$\\
 \hline
41 & $n-1$ & $n-1$ & $n-1$ & 191 & $n-1$ & $n-1$ & $n-1$\\
 \hline
43 & $\lceil\frac{n}{2}\rceil+1$ & $\lceil\frac{n}{2}\rceil$ & $\lceil\frac{n}{2}\rceil-1$ & 193 & $\lceil\frac{n}{2}\rceil-1$ & $\lceil\frac{n}{2}\rceil-1$ & $\lceil\frac{n}{2}\rceil-1$ \\
 \hline
47 & $n+1$ & $n$ & $n-1$ & 197 & $n$ & $n$ & $n-1$\\
 \hline
53 & $n$ & $n$ & $n-1$ & 199 & $\lceil\frac{n}{2}\rceil$ & $\lceil\frac{n}{2}\rceil-1$ & $\lceil\frac{n}{2}\rceil-1$ \\
 \hline
59 & $n$ & $n-1$ & $n-1$ & 211 & $n-1$ & $n-1$ & $n-1$\\
 \hline
61 & $n-1$ & $n-1$ & $n-1$ & 223 & $n+1$ & $n$ & $n-1$\\
 \hline
67 & $n+1$ & $n$ & $n-1$ & 227 & $n+1$ & $n$ & $n-1$\\
 \hline
71 & $n-1$ & $n-1$ & $n-1$ & 229 & $n$ & $n-1$ & $n-1$\\
 \hline
73 & $n$ & $n$ & $n-1$ & 233 & $n$ & $n$ & $n-1$\\
 \hline
79 & $n$ & $n-1$ & $n-1$ & 239 & $n$ & $n-1$ & $n-1$\\
 \hline
83 & $n+1$ & $n$ & $n-1$ & 241 & $n-1$ & $n-1$ & $n-1$\\
 \hline
89 & $n$ & $n-1$ & $n-1$ & 251 & $\lfloor\frac{n}{2}\rfloor$ & $\lceil\frac{n}{3}\rceil-1$ & $\lceil\frac{n}{3}\rceil-1$\\
 \hline
97 & $n$ & $n$ & $n-1$ & 257 & $\lceil\frac{n}{2}\rceil$ & $\lceil\frac{n}{2}\rceil$ & $\lceil\frac{n}{2}\rceil-1$\\
 \hline
101 & $\lceil\frac{n}{2}\rceil-1$ & $\lceil\frac{n}{2}\rceil-1$ & $\lceil\frac{n}{2}\rceil-1$ & 263 & $n+1$ & $n$ & $n-1$\\
 \hline
103 & $n+1$ & $n$ & $n-1$ & 269 & $n$ & $n-1$ & $n-1$\\
 \hline
107 & $\lceil\frac{n}{2}\rceil+1$ & $\lceil\frac{n}{2}\rceil$ & $\lceil\frac{n}{2}\rceil-1$ & 271 & $n-1$ & $n-1$ & $n-1$\\ 
 \hline
109 & $n$ & $n-1$ & $n-1$ & 277 & $n$ & $n$ & $n-1$\\
 \hline
113 & $n$ & $n$ & $n-1$ & 281 & $n-1$ & $n-1$ & $n-1$\\
 \hline
\end{tabular}
\caption{Values of $f_p(x,0,n)$ by varying $p$ and $x$}
\label{tab:4.1}
\end{table} 
Where in the last column $x\geq2$. \\
\subsection{When $q\equiv 0\mod 2$}
When $q\equiv 0\mod 2$ and $x<2$ it could be the most anomalous case. Because as we have seen in my last papers, the repeating last digits of a tetration which base is a power of 2 has a lot of sub-cases. So it may be complicated to find a generic "nice" formula when $x<2$.
\section{Conclusions}
We are very near to a proof for a formula which finds the minimum hyper-exponent $u$ of a tetration with a generic base $q$ such that the last $n$ digits of the tetration after the $u$-th one are the same.


\begin{thebibliography}{1}
\bibitem{tet}Womack David, Repeated powers: The operation of tetration \textit{Mathematics in School} volume 42, number 4, pages 38--40, JSTOR, 2013.
\bibitem{tdn}Salvatore Damantino, Teoria dei numeri, \textit{Teoria dei numeri}, pages 116--125, scienza express, 2018.
\end{thebibliography}
\end{document}